\documentclass[reqno,twoside,a4paper,10pt]{amsart}
\usepackage[latin1]{inputenc}
\usepackage[english]{babel}
\usepackage{latexsym,amsthm,amsmath,amscd,amsxtra,amssymb,amsxtra,exscale,mathrsfs,bbm,ifthen}
\usepackage[colorlinks,linkcolor={blue},citecolor={blue},urlcolor={red},]{hyperref}

\theoremstyle{plain}
\newtheorem{theorem}{Theorem}[section]

\newtheorem{lemma}[theorem]{Lemma}
\newtheorem{proposition}[theorem]{Proposition}

\theoremstyle{definition}

\theoremstyle{remark}
\newtheorem{remark}[theorem]{Remark}

\newtheorem*{remark*}{Remark}

\newtheorem{example}[theorem]{Example}

\newcommand{\NN}{{\mathbb N}}
\newcommand{\ZZ}{{\mathbb Z}}

\newcommand{\RR}{{\mathbb R}}
\newcommand{\CC}{{\mathbb C}}
\newcommand{\KK}{{\mathbb K}}

\newcommand{\EE}{{\mathbb E}}
\newcommand{\PP}{{\mathbb P}}
\newcommand{\cF}{{\mathcal F}}

\newcommand{\cH}{{\mathcal H}}

\newcommand{\al}{\alpha}
\newcommand{\ga}{\gamma}
\newcommand{\eps}{\varepsilon}
\newcommand{\la}{\lambda}
\newcommand{\om}{\omega}
\newcommand{\Om}{\Omega}

\newcommand{\embed}{\hookrightarrow}

\newcommand{\cB}{{\mathcal B}}
\newcommand{\norm}[1] {\| #1 \|}
\newcommand{\lrnorm}[1]{\left\| #1 \right\|}

\newcommand{\biggnorm}[1]{\biggl\| #1 \biggr\|}
\newcommand{\sfrac}[2] { {{}^{#1}\!\!/\!{}_{#2}}}
\newcommand{\einhalb} {\sfrac{1}{2}}

\newcommand{\as}{{\ga_\infty}}

\renewcommand{\Re}[1] {{\rm Re}#1}
\renewcommand{\Im}[1] {{\rm Im}#1}
\newcommand{\suchthat}{\,\colon\;}
\newcommand{\SCP}[2]{ {\rm(SCP)}_{(#1,#2)}}
\renewcommand{\span}{{\rm span}}
\newcommand{\eins}{{\mathbbm 1}}
\newcommand{\RB}[2][\relax]{ {
   \mathcal{R}\left(\{ #2 \}\right)_{#1}}}

\newcommand{\essinf}{{\rm ess}\,\inf}
\newcommand{\esssup}{{\rm ess}\,\sup}

\newcommand {\comment}[1]{\relax}

\begin{document}

\title{A stochastic Datko-Pazy theorem}

\author{Bernhard Haak \and Jan van Neerven \and Mark Veraar
}
\date\today
\thanks{The authors gratefully acknowledge financial support by
          a `VIDI subsidie' (639.032.201) in the `Vernieuwingsimpuls'
          programme of the Netherlands Organization for Scientific
          Research (NWO). The second named author is also supported by
          a Research Training Network (HPRN-CT-2002-00281).}
\maketitle

\def\leftmark{\hfill{\sc A stochastic Datko-Pazy theorem}\hfill}
\def\rightmark{\hfill{\sc B.H. Haak, J.M.A.M. van Neerven, M.C. Veraar}\hfill}

\newcommand{\contact}{\noindent Delft Institute of Applied Mathematics\\
 Technical University of Delft \\ P.O. Box 5031\\ 2600 GA Delft\\The
 Netherlands \\[2ex]
        {\it E-mail addresses:} {\tt B.H.Haak@math.tudelft.nl}\\
\phantom{\it E-mail addresses:} {\tt J.vanNeerven@math.tudelft.nl}\\
\phantom{\it E-mail addresses:} {\tt M.C.Veraar@math.tudelft.nl}
}

\def\keywords{semigroups, Datko-Pazy theorem, stochastic Cauchy 
   problem, invariant measures, perturbation theory Riesz sequences, 
   almost summing operators, $\ga$--radonifying operators }
\def\subjclass{Primary:  47D06, 28C20, Secondary: 46B09, 46B15, 47N30}

\begin{abstract}
\noindent Let $H$ be a Hilbert space and $E$ a Banach space. In this note we
present a sufficient condition for an operator $R: H\to E$ to be
$\ga$--radonifying in terms of Riesz sequences in $H$. This result is applied
to recover a result of Lutz Weis and the second named author on the
$R$-boundedness of resolvents, which is used to obtain a Datko-Pazy type
theorem for the stochastic Cauchy problem. We also present some
perturbation results.

\bigskip

\noindent {\small {\bf Subject classifications:} 
           \hspace*{\labelsep} \subjclass}\\[1ex]
\noindent {\small {\bf Key words:}\hspace*{\labelsep} \keywords}
\end{abstract}

\section{Introduction}

The well-known Datko-Pazy theorem
states that if $(T(t))_{t\ge 0}$ is a strongly
continuous semigroup on a Banach space $E$ such that all orbits
$T(\cdot)x$ belong to the space $L^p(\RR_+, E)$ for some $p\in
[1,\infty)$, then $(T(t))_{t\ge 0}$ is uniformly exponentially
stable, or equivalently, there exists an $\eps >0$ such that
all orbits $t \mapsto e^{\eps t} T(t)x$ belong to $L^p(\RR_+, E)$.
For $p=2$ and Hilbert spaces $E$ this result is due to Datko \cite{Datko},
and the general case was obtained by Pazy \cite{Pazy:Datko}.

In this note we prove a stochastic version of the Datko-Pazy
theorem for spaces of
$\ga$--radonifying operators (cf. Section~\ref{sec:2}).
Let us denote by $\ga(\RR_+,E)$ the space of all strongly measurable
functions $\phi: \RR_+ \to E$ for which the integral operator
\[
f \mapsto \int_0^\infty f(t) \phi(t)\,dt
\]
is well-defined and $\ga$-radonifying from
$L^2(\RR_+)$ to $E$.

\let\ALTERWERT\thetheorem
\def\thetheorem{\arabic{section}.\arabic{theorem}a}
\begin{theorem}[Stochastic Datko-Pazy Theorem, first version]\label{thm:DatkoPazy}
Let $A$ be the generator of a strongly continuous semigroup
$(T(t))_{t\ge  0}$ on a Banach space $E$. The following assertions are equivalent:
\begin{enumerate}
\item For all $x\in E$, $T(\cdot) x\in \ga(\RR_+, E)$.
\item There exists an $\eps>0$ such that for all $x\in E$, $t\mapsto e^{\eps t} T(t) x\in
  \ga(\RR_+, E)$.
\end{enumerate}
\end{theorem}
\let\thetheorem\ALTERWERT
If $E$ is a Hilbert space, $\ga(\RR_+, E) = L^2(\RR_+, E)$ and
Theorem~\ref{thm:DatkoPazy} is equivalent to the Datko's theorem mentioned above.

\medskip

As explained in \cite{vanNeervenWeis:asymptotic-scp},
$\ga$--radonifying operators play an important role in the study of the
following stochastic abstract Cauchy problem on $E$:
\[
\SCP{A}{B} \qquad \left\{
  \begin{array}{lcl}
    dU(t) &=& A\, U(t)\, dt + B \,dW_H(t), \qquad t\ge 0,\\
    U(0)  &=& 0.
  \end{array}
 \right.
\]
Here, $H$ is a separable Hilbert space, $B \in \cB(H, E)$ is a bounded operator, and $W_H$ is an $H$-cylindrical Brownian motion.

Theorem~\ref{thm:DatkoPazy} can be reformulated in terms of invariant measures for $\SCP{A}{B}$ as follows.

\let\ALTERWERT\thetheorem
\def\thetheorem{\arabic{section}.\arabic{theorem}b}
\setcounter{theorem}{0}
\begin{theorem}[Stochastic Datko-Pazy theorem, second version]
\label{thm:DatkoPazy2}
With the above notations, the following assertions are equivalent:
\begin{enumerate}
\item For all rank one operators $B\in \cB(H, E)$,
     the problem $\SCP{A}{B}$ admits an invariant measure.
\item There exists an $\eps>0$ such that for all rank one operators \ $B\in \cB(H, E)$,\
     the problem \ $\SCP{A+\eps}{B}$ admits an invariant measure.

\end{enumerate}

\end{theorem}
\let\thetheorem\ALTERWERT

For unexplained terminology
and more information on the stochasic Cauchy problem and
invariant measures we refer to
\cite{DaPratoZabczyk,vanNeervenWeis:stoch-int,vanNeervenWeis:asymptotic-scp}.

\section{Riesz bases and $\ga$-radonifying operators}\label{sec:2}

Let $\cH$ be a Hilbert space and $E$ a Banach space.
Let $(\ga_n)_{n\ge 1}$ be a sequence of independent standard Gaussian
random variables on a probability space $(\Om, \cF, \PP)$.
A bounded linear operator $R: \cH \to E$ is called {\em almost summing} if
\[
\norm{R}_{\as(\cH, E)} := \sup \biggnorm{ \sum_{n=1}^N \ga_n R h_n }_{L^2(\Omega, E)}  < \infty,
\]
where the supremum is taken over all $N\in \NN$ and all orthonormal
systems $\{ h_1, \ldots, h_N\}$ in $\cH$. Endowed with this
norm, the space $\as(\cH, E)$ of all almost summing operators
is a Banach space. Moreover, $\as(\cH,E)$ is an operator ideal in $\cB(\cH,E)$.
The closure of the finite rank operators in $\as(\cH,E)$ will be denoted by
$\ga(\cH,E)$. Operators belonging to this space are called {\em $\ga$-radonifying}. Again
$\ga(\cH,E)$ is an operator ideal in $\cB(\cH,E)$.

\comment{
\begin{proof}
Let $(R_n)_{n\geq 1}$ be a Cauchy sequence in $\as(\cH, E)$ Consequently,
there exists a constant $M>0$ such that $\norm{R_n}_{\as(\cH,E)} \le M$
for all $n\in \NN$. Since $R_n$ is also Cauchy in $\cB(\cH, E)$, it
converges towards a bounded operator $R\in \cB(\cH,E)$. Let $\eps >0$. Then
there is a $N_0$ such that for $M,N\ge N_0$ we have
$\norm{R_N-R_M}_{\as(\cH,E)} < \eps$. Let $(h_n)_{n=1}^m$ be a finite
orthonormal system. Then by Fatou's lemma,
\[
     \biggnorm{ \sum_{n=1}^m \ga_n (R-R_N) h_n }_{L^2(\Omega, E)}
\leq \liminf_{M\to \infty}\biggnorm{ \sum_{n=1}^m \ga_n (R_M-R_N)
      h_n}_{L^2(\Omega, E)}<\eps.
\]
This being valid for all finite orthonormal systems, we obtain
$R\in \as(\cH,E)$ and $\norm{R-R_N}_{\as(\cH,E)} < \eps$ for all $N>N_0$, implying the
completeness of $\as(\cH,E)$.
\end{proof}
}

Let us now assume that $\cH$ is a separable Hilbert space. Under this assumption one has $R\in \as(\cH,E)$ if and only if for some (every) orthonormal basis $(h_n)_{n\geq 1}$ for $\cH$,
\[
M := \sup_{N\geq 1}\biggnorm{\sum_{n=1}^N \ga_n Rh_n}_{L^2(\Om,E)} < \infty.
\]
In that case, $\norm{R}_{\as(\cH,E)} = M$.
Furthermore, one has $R\in
\ga(\cH,E)$ if and only if for some (every) orthonormal basis
$(h_n)_{n\geq 1}$ for $\cH$, $\sum_{n\geq 1} \ga_n Rh_n$ converges in
$L^2(\Om, E)$. In that case,
\[
\norm{R}_{\ga(\cH,E)} = \biggnorm{\sum_{n\geq 1} \ga_n Rh_n }_{L^2(\Om,E)}.
\]
If $E$ does not contain a closed subspace isomorphic to $c_0$, then by a
result of Hoffmann-J\o{}rgensen and Kwapie\'n (cf. \cite[Theorem 9.29]{LT}),
$\ga(\cH, E) = \as(\cH, E)$.

We will apply the above notions to the space $\cH = L^2(\RR_+,H)$
where $H$ is a separable Hilbert space. For an operator-valued function
$\phi: \RR_+ \to \cB(H,E)$ which is {\em $H$-strongly measurable} in the sense that $t\mapsto \phi(t)h$ is strongly measurable for all $h\in H$, and  {\em weakly square integrable} in the sense that $t\mapsto \phi^*(t)x^*$ is square Bochner integrable for all $x^*\in E^*$, let
$R_\phi \in \cB(L^2(\RR_+, H), E)$ be defined as the Pettis integral
operator  $$R_\phi(f) := \int_{\RR_+}\phi(t) f(t)\,dt.$$
We say that $\phi \in \ga(\RR_+, H, E)$ if $R_\phi \in \ga( L^2(\RR_+,
H), E)$ and write $$\norm{\phi}_{\ga(\RR_+, H,E)} := \norm{R_\phi}_{\ga( L^2(\RR_+,H), E)}.$$
If $H=\KK$, where $\KK =\RR$ or $\CC$ is the underlying
scalar field, we write $\ga(\RR_+, E)$ for $\ga(\RR_+, H, E)$. For
almost summing operators we use an analogous notation.

For more information we refer to
\cite{DiestelJarchowTonge,KaltonWeis:square-function-est,vanNeervenWeis:stoch-int,vanNeervenWeis:asymptotic-scp}.

\paragraph{Hilbert and Bessel sequences.}
Let $\cH$ be a Hilbert space and
$I\subseteq \ZZ$ an index set. A sequence $(h_i)_{i\in I}$ in $\cH$ is said
to be a {\em Hilbert sequence} if there exists a constant $C>0$ such that for
all scalars $(\al_i)_{i\in I}$,
\[
       \biggl(\biggnorm{ \sum_{i\in I} \al_i h_i }^2\biggr)^\einhalb
\leq C \biggl(\sum_{i\in I} |\al_i|^2\biggr)^\einhalb.
\]
The infimum of all admissible constants $C>0$ will be denoted by $C_H(\{h_i
\suchthat i\in I\})$. A Hilbert sequence that is a Schauder basis
is called a {\em Hilbert basis} (cf. \cite[Section 1.8]{Young:nonharmonic}).

The sequence $(h_i)_{i\in I}$ is said to be a {\em Bessel sequence} if there
exists a constant $c>0$ such that for all scalars $(\al_i)_{i\in I}$,
\[
       c \biggl(\sum_{i\in I}|\al_i|^2 \biggr)^\einhalb
  \leq \biggl(\biggnorm{ \sum_{i\in I} \al_i h_i }^2\biggr)^\einhalb.
\]
The supremum of all admissible constants $c>0$ will be denoted by
$C_B(\{h_i\suchthat i\in I\})$. Notice that every Bessel sequence is
linearly independent. A Bessel sequence that is a
Schauder basis is called a {\em Bessel basis}. A sequence
$(h_i)_{i\in I}$ that is a Bessel sequence and a Hilbert sequence
is said to be a {\em Riesz sequence}. A sequence $(h_i)_{i\in I}$ that is a
Bessel basis and a Hilbert basis is said to be a {\em Riesz basis} (cf.
\cite[Section 1.8]{Young:nonharmonic}).

In the above situation if it is clear which sequence in $\cH$ we refer to, we
use the short-hand notation $C_H$ and $C_B$ for
$C_H(\{h_i\suchthat i\in I\})$ and $C_B(\{h_i\suchthat i\in
I\})$.

In the next results we study the relation between $\ga$--radonifying
operators and Hilbert and Bessel sequences.

\begin{proposition}\label{prop:Hilbert-sequence-estimate}

Let $(f_n)_{n\geq
1}$ be a Hilbert sequence in $\cH$.
\begin{enumerate}
\item \label{item:Hilbert-sequence-estimate-a}
If $R\in \as(\cH,E)$, then
\begin{equation}\label{eq:HS-almost-summing}
\sup_{N\geq 1}\biggnorm{ \sum_{n =1}^N \ga_n R f_n }_{L^2(\Om, E)}
  \leq C_H \; \norm{R}_{\as(\cH,E)}.
\end{equation}
\item \label{item:Hilbert-sequence-estimate-b}
If $R\in \ga(\cH,E)$, then $\sum\limits_{n\geq 1} \ga_n R f_n$ converges in
$L^2(\Om, E)$ and
\begin{equation}\label{eq:HS-gamma}
\biggnorm{ \sum_{n \geq 1} \ga_n R f_n }_{L^2(\Om, E)} \leq C_H \;
\norm{R}_{\ga(\cH,E)}.
\end{equation}
\end{enumerate}
\end{proposition}
\begin{proof}
\ref{item:Hilbert-sequence-estimate-a}: Fix $N\ge 1$ and let $\{h_1, \ldots, h_N\}$ be an
orthonormal system in $\cH$. Since $(f_n)_{n\ge 1}$ is a Hilbert sequence there
is a unique $T\in \cB(\cH)$ such that $T h_n = f_n$ for $n=1,\ldots,N$
and $Tx=0$ for all $x\in \{h_1,\ldots,h_N\}^\perp$. Moreover,
$\norm{T} \le C_H$.

\comment{
Indeed, define a linear operator $T: \span(h_n, n\geq 1)\to \cH$ by
$T h_n = f_n$.  For $h = \sum_{n=1}^N a_n h_n$ we have
\[
  \norm{T h}
= \biggnorm{\sum_{n=1}^N a_n f_n  }\leq C_H \biggl(\sum_{n=1}^N
    |a_n|^2\biggr)^{\einhalb}
= C_H \norm{h}.
\]
It follows that $T$ has a continuous
extension to an operator $T\in \cB(\cH)$ of norm $C_H$.
}

By the right ideal property we have $R\circ T\in \as(\cH,E)$ and, for all $N\ge 1$,
\[
\biggnorm{ \sum_{n =1}^N \ga_n R f_n }_{L^2(\Om, E)} = \biggnorm{ \sum_{n =1}^N
\ga_n R T h_n }_{L^2(\Om, E)}\leq \norm{R\circ T}_{\as(\cH,E)}\leq C_H \;
\norm{R}_{\as(\cH,E)}.
\]
\ref{item:Hilbert-sequence-estimate-b}: This is proved in a similar way.

\end{proof}

\begin{proposition}\label{prop:Bessel-sequence-estimate}

Let
$(f_n)_{n\geq 1}$ be a Bessel sequence in $\cH$ and let
$\cH_f$ denote its closed linear span.
\begin{enumerate}
\item\label{item:Bessel-basis-estimate-a}
If $\sup\limits_{N\ge 1}\lrnorm{ \sum\limits_{n =1}^N \ga_n R f_n }_{L^2(\Om, E)}<\infty$,
then $R\in \as(\cH_f,E)$ and
\begin{equation}\label{eq:BB-almost-summing}
     \norm{R}_{\as(\cH_f, E)}
\leq C_B^{-1} \; \sup_{N\ge1} \biggnorm{ \sum_{n =1}^N
      \ga_n R f_n }_{L^2(\Om, E)}.
\end{equation}
\item\label{item:Bessel-basis-estimate-b}
If $\sum\limits_{n\geq 1} \ga_n R f_n$ converges in $L^2(\Om, E)$, then
$R\in \ga(\cH_f,E)$ and
\begin{equation}\label{eq:BB-gamma}
     \norm{R}_{\ga(\cH_f,E)}
\leq C_B^{-1} \; \biggnorm{ \sum_{n \geq 1} \ga_n R
       f_n}_{L^2(\Om, E)}.
\end{equation}
\end{enumerate}
\end{proposition}
\begin{proof}
Let $(h_n)_{n\ge 1}$ an orthonormal basis for $\cH_f$. Since $(f_n)_{n\geq 1}$ is a Bessel sequence there
is a unique $T\in \cB(\cH,E)$ such that $T f_n = h_n$ and
$Tx = 0$ for $x \in\cH_f^\perp$. Notice that
$\norm{T} \le C_B^{-1}$.
On the linear span $\cH_0$ of the sequence
$(f_n)_{n\ge 1}$ we define an inner product by
$[x,y]_T := [Tx, Ty]_\cH$. Note that this is well defined
by the linear independence of the sequence $(f_n)_{n\ge 1}$.
Let $\cH_T$ denote the Hilbert space completion of $\cH_0$ with respect to
$[\cdot,\cdot]_T$. The identity mapping on $\cH_f$ extends
to a bounded operator $j:\cH_f\embed \cH_T$ with norm $\norm{ j } \le C_B^{-1}$.
Clearly, $(j f_n)_{n\geq 1}$ is an orthonormal sequence in $\cH_T$ with dense
span, and therefore it is an
orthonormal basis for $\cH_T$. It is elementary to verify that
the assumption on $R$ may now be translated as saying that $R$ extends
in a unique way to an almost summing operator (in part
\ref{item:Bessel-basis-estimate-a}), respectively a $\ga$-radonifying
operator (in part \ref{item:Bessel-basis-estimate-b}), denoted by $R_T$, from $\cH_T$ to $E$.
We estimate
\begin{eqnarray*}
    \biggnorm{\sum_{n\ge 1} \al_n jh_n }_{\cH_T}
= \biggnorm{\sum_{n\ge 1} \al_n T h_n }_{\cH}
\; \le \; C_B^{-1} \; \biggnorm{\sum_{n\ge 1} \al_n h_n }_{\cH}
= C_B^{-1} \; \biggl(\sum_{n\ge 1} |\al_n|^2 \biggr)^\einhalb.
\end{eqnarray*}
From this we deduce that $(jh_n)_{n\ge 1}$ is a Hilbert sequence in $\cH_T$ with constant
$\le C_B^{-1}$.
Hence we may apply Proposition~\ref{prop:Hilbert-sequence-estimate} to the operator
$R_T: \cH_T\to E$ and the Hilbert sequence $(jh_n)_{n\ge 1}$ in $\cH_T$ to
obtain the result.
\end{proof}

As a consequence of the above results we obtain:

\begin{theorem}\label{thm:zusammenfassung}

Let
$(f_n)_{n\geq 1}$ be a Riesz basis in the Hilbert space $\cH$.
\begin{enumerate}
\item
One has $R\in \as(\cH,E)$ if and only if $\sup\limits_{N\geq 1}\biggnorm{ \sum\limits_{n=1}^N \ga_n R
f_n }_{L^2(\Om, E)}\!\!\!<\infty$. In that case \eqref{eq:HS-almost-summing} and
 \eqref{eq:BB-almost-summing} hold.
\item One has $R\in \ga(\cH,E)$ if and only if $\sum\limits_{n\geq 1} \ga_n R f_n$
  converges in $L^2(\Om, E)$. In that case \eqref{eq:HS-gamma} and
  \eqref{eq:BB-gamma} hold.
\end{enumerate}
\end{theorem}

The following well-known lemma identifies a class of Riesz sequences
in $L^2(\RR)$. For convenience we include the short proof from
\cite[Theorem 2.1]{CasazzaChristensenKalton}. Let $\mathbb{T}$ be the
unit circle in $\mathbb{C}$.

\begin{lemma}\label{lem:CCK}

Let $f\in L^2(\RR)$ and define the sequence $(f_n)_{n\in \mathbb{Z}}$ in
$L^2(\RR)$ by $f_n(t) = e^{2\pi  n i t} f(t)$. Define $F:\mathbb{T}\to
\overline{\RR}$ as
\[
F(e^{2\pi i t}) := \sum_{k\in \mathbb{Z}} |f(t+ k)|^2\biggr.
\]
\begin{enumerate}
\item The sequence $(f_n)_{n\in \mathbb{Z}}$ is a Bessel sequence in $L^2(\RR)$ if and only if
there exists a constant $A>0$ such that $A\leq F(e^{2\pi i t} )$ for
almost all $t\in [0,1]$.
\item The sequence  $(f_n)_{n\in \mathbb{Z}}$ is a Hilbert sequence in $L^2(\RR)$ if and only if
there exists a constant $B>0$ such that $F(e^{2\pi i t} ) \leq B$ for
almost all $t\in [0,1]$.
\end{enumerate}
In these cases, $C_B^{2} = \essinf F$ and $C_H^2= \esssup F$ respectively.
\end{lemma}
\begin{proof}
Both assertions are obtained by observing that for $I\subseteq\ZZ$ and
$(a_n)_{n\in I}$ in $\mathbb{C}$ we may write
\begin{eqnarray*}
    \biggnorm{\sum_{n\in I} a_n f_n}_{L^2(\RR)}^2
&=& \sum_{k\in \ZZ} \int_{k }^{(k+1)} \biggl| \sum_{n\in I} a_n
      e^{2\pi n  i t} f(t)\biggr|^2 \, dt
\\ &=&  \sum_{k\in \ZZ} \int_0^{1} \biggl| \sum_{n\in I} a_n e^{2\pi n i t}
      f(t+k)\biggr|^2 \, dt
= \int_0^{1} \biggl| \sum_{n\in I} a_n e^{2\pi n i t}\biggr|^2
    F(e^{2\pi i t}) \, dt.
\end{eqnarray*}
\end{proof}

The following application of Lemma~\ref{lem:CCK} will be used below.

\begin{example}\label{ex:hilbert-sequence-strip-type2}

Let $\rho\in [0,1)$ and $a>0$. For $n\in \ZZ$ let
\[
      f_n(t) = e^{-a t+ 2\pi (n+\rho) i t} \eins_{[0,\infty)}(t).
\]
Then $(f_n)_{n \in \ZZ}$ is a Riesz sequence in $L^2(\RR)$ with constants
$C_B^2 = \tfrac{e^{-2a}}{e^{2a}-1}$ and $C_H^2 =
\tfrac{e^{2a}}{e^{2a}-1}$.
Indeed, let $f(t) := e^{-a t+ 2\pi \rho i t} \eins_{[0,\infty)}(t)$. For all $t\in [0,1)$,
\[
   F(e^{2\pi i t}) = \sum_{k\in \mathbb{Z}} |f(t+ k)|^2
 = \sum_{k=0}^\infty e^{-2a(t+k)}
 = \frac{e^{2a(1-t)}}{e^{2a}-1}.
\]
Now Lemma~\ref{lem:CCK} implies the result.
\end{example}

\begin{remark}
Necessary and sufficient conditions on the complex coefficients $c_n$ and
$\la_n$ with $\Re\la_n >0$ in order that the functions $z\mapsto c_n \exp(-\la_n z)$
form a Riesz sequence can be found in \cite[Section 10.3]{NiPa} and \cite{JZ}.
\end{remark}

\section{Main results}

In this section we use Proposition~\ref{prop:Hilbert-sequence-estimate} to obtain an
alternative proof of  \cite[Theorem 3.4]{vanNeervenWeis:asymptotic-scp}
on the $R$--boundedness of certain Laplace transforms. This result is
applied to strongly continuous semigroups to obtain estimates for the
abscissa of $R$--boundedness of the resolvent. From this we deduce
Theorem~\ref{thm:DatkoPazy} as well as  bounded perturbation results for the existence of solutions and invariant measures
for the problem $\SCP{A}{B}$.

Let $(r_n)_{n\ge 1}$ be a Rademacher sequence on a probability space
$(\Om,{\mathscr F},\PP)$.
A family of operators ${\mathscr T}\subseteq \cB(E)$ is called
{\em $R$-bounded}
if there exists a constant $C>0$ such that for all $N\ge 1$ and all sequences
$(T_n)_{n=1}^N\subseteq {\mathcal T}$ and
$(x_n)_{n=1}^N\subseteq E$ we have
$$ \EE\Bigl\Vert \sum_{n=1}^N r_n T_n x_n\Bigr\Vert^2 \le C^2
 \EE\Bigl\Vert \sum_{n=1}^N r_n x_n\Bigr\Vert^2.
$$
The least possible constant $C$ is called the
{\em $R$-bound} of ${\mathscr T}$, notation ${\mathscr R}({\mathscr T})$.
Clearly, every $R$-bounded family ${\mathscr T}$ is uniformly bounded and
$\sup_{T\in {\mathscr T}} \Vert T\Vert \le {\mathscr R}({\mathscr T})$.

Following \cite{vanNeervenWeis:asymptotic-scp}, for an operator $T\in
\cB(L^2(\RR_+), E)$ we define the {\em Laplace transform}
$\widehat{T}:\{\la\in \CC \suchthat \Re{\la} >0\} \to E$ as
\[
   \widehat{T}(\la) := T e_{\la}.
\]
Here $e_{\la}\in L^2(\RR_+)$ is given by $e_{\la}(t) =
e^{-\la t}$. For a Banach space $F$ and a bounded operator
$\Theta: F\to \cB(L^2(\RR_+), E)$ we define the
{\em Laplace transform} $\widehat{\Theta}: \{ \la\in \CC \suchthat
\Re{\la}>0\}\to \cB(F,E)$ as
\[
   \widehat{\Theta}(\la) y
:=  \widehat{\Theta y}(\la) \qquad  \Re{\la}>0, \ y\in F.
\]

The following result is a slight refinement of \cite[Theorem 3.4]{vanNeervenWeis:asymptotic-scp}. The main novelty is the simple proof
of the estimate \eqref{eq:estimate}.

\begin{theorem}\label{thm:RboundednessLaplace}

Let $F$ be a Banach space. Let $\Theta:F\to \as(L^2(\RR_+), E)$ be a bounded
operator and let $\delta>0$. Then $\widehat{\Theta}$ is $R$--bounded on the
half-plane $\{ \la\in \CC \suchthat \Re{\la} > \delta \}$ and there
exists a universal constant $C$ such that
\[
     \RB{ \widehat{\Theta}(\la)\suchthat  \Re{\la} \geq \delta}
\leq \norm{\Theta} \frac{C}{\sqrt{\delta}}.
\]
\end{theorem}
\begin{proof}
Let $\delta >0$. 
Consider the set $\{\la\in\CC: \  \Re{\la} =\delta\}$. Fix $\sigma\in [\sfrac{\delta}2,
\sfrac{3}2\delta]$

and $\rho\in [0,1)$. For $n\in \ZZ$ let $g_n:\RR_+\to \CC$ be given by
\[
  g_n(t) = e^{-\sigma t+ (n+\rho) \delta  i t}.
\]

By substitution, this reduces to
Example~\ref{ex:hilbert-sequence-strip-type2}, whence $(g_n)_{n\geq 1}$ is a
Riesz sequence in $L^2(\RR_+)$ with constant $0<C_H \leq
\bigl(\tfrac{C}{\delta}\bigr)^\einhalb$ where $C := 2\pi
\tfrac{e^{2\pi}}{e^{2\pi}-1}$. For $y\in F$, we may apply
Proposition~\ref{prop:Hilbert-sequence-estimate} to obtain
\begin{equation}
\begin{aligned}
\label{eq:estimate}
      \biggnorm{\sum_{n=-N}^N \ga_n \widehat{\Theta}(\sigma -
        (n+\rho)\delta i) y }_{L^2(\Om,E)}
& =  \biggnorm{\sum_{n=-N}^N \ga_n (\Theta y) g_n }_{L^2(\Om,E)}\\
&\leq C_H \norm{\Theta y}_{\as(\Om, E)}
 \leq \Bigl(\frac{C}{\delta}\Bigr)^\einhalb \norm{\Theta} \, \norm{y}.
\end{aligned}
\end{equation}
The rest of the proof follows the lines in
\cite{vanNeervenWeis:asymptotic-scp}.

\end{proof}

In what follows we let $(T(t))_{t\ge 0}$ be a strongly continuous semigroup on $E$ with
generator $A$. We recall from
\cite{vanNeervenWeis:stoch-int, vanNeervenWeis:asymptotic-scp}
that the problem $\SCP{A}{B}$ admits a (unique) solution if and only
if $T(\cdot)B$ belongs to $\ga([0,T],H,E)$ for some (all) $T>0$.
Furthermore, an invariant measure exists if and only if
$T(\cdot)B$ belongs to $\ga(\RR_+,H,E)$.

The next theorem improves \cite[Theorem 1.3]{vanNeervenWeis:asymptotic-scp},
where the bound $s_R(A)\le 0$ was obtained.

\begin{theorem}\label{thm:gammaDatkoPazyResolvent}

Assume that for all $x\in E$, $T(\cdot)x\in \as(\RR_+,E)$. Then
$s_R(A) <0$, i.e., there exists an $\eps>0$ such that $\{R(\la,
A)\suchthat  \Re{\la} \geq -\eps\}$ is $R$--bounded.
\end{theorem}
\begin{proof}
By the closed graph theorem there exists an $M>0$ such that $\norm{
  T(\cdot)x }_{\as(\RR_+, E)} \le M \norm{x}$. By
Theorem~\ref{thm:RboundednessLaplace}, $\{\la \in\CC: \ \Re\la >0\}\subseteq \varrho(A)$ and
\begin{equation}\label{eq:abfall-der-R-schranke}
     \RB{ R(\la, A)\suchthat  \Re{\la} \geq \delta }
\leq \frac{c}{\sqrt{\delta}}
\end{equation}
for all $\delta>0$, where $c := CM$ with $C$ the universal constant
of Theorem~\ref{thm:RboundednessLaplace}.
The following standard argument shows that this implies the bound
\begin{equation}\label{sA}
s(A) \leq -\frac{1}{4c^2}.
\end{equation}
Choose $\delta>0$ and let $\mu\in \sigma(A)$ be such that
$\Re{\mu} > s(A)-\delta$.  With $\la = \frac{1}{4 c^2} + i\,
\Im{\mu}$ it follows that
\[
     \frac{1}{4 c^2} - s(A) +\delta
\geq  {\rm dist}(\la, \sigma(A))
\geq \frac{1}{\norm{R(\la, A)}}
\geq \frac{ \sqrt{\Re{\la} } }{c}  =  \frac{1}{2 c^2}.
\]
Thus $s(A) \leq -\frac{1}{4c^2}+\delta$. Since $\delta>0$ was
arbitrary, this gives \eqref{sA}.

Now let $\eps_0 := \frac{1}{4c^2}$. For $\la$ with $-\eps_0<
\Re{\la}<3\eps_0$ we may write
\[
R(\la,A)  = \sum_{n\geq 0} (\eps_0 -  \Re{ \la})^n R(\eps_0 + i
\Im{\la}, A)^{n+1}.
\]
Fix $0<\eps<\eps_0$. We claim that
$\{ R(\la,A)\suchthat  \Re{\la} = -\eps\}$ is
$R$--bounded. To see this let $(r_k)_{k=1}^K$ be a Rademacher sequence
on $(\Om,{\mathscr F},\PP)$, let $(\la_k)_{k=1}^K$ be such that
$\Re{ \la_k } = -\eps$, and let $(x_k)_{k=1}^K$ be a sequence in $E$. We may
estimate
\begin{eqnarray*}
   \biggnorm{\sum_{k=1}^K r_k R(\la_k, A) x_k }_{L^2(\Om,E)}
& = & \biggnorm{\sum_{n\geq 0} \sum_{k=1}^K r_k  (\eps_0 + \eps)^n
        R(\eps_0 + i \Im{\la_k}, A)^{n+1} x_k }_{L^2(\Om,E)}\\
&\leq&\sum_{n\geq 0} (\eps_0 + \eps)^n \biggnorm{\sum_{k=1}^K r_k
        R(\eps_0 + i \Im{\la_k}, A)^{n+1} x_k }_{L^2(\Om,E)}\\
&\leq& \sum_{n\geq 0} (\eps_0 + \eps)^n \biggl( \frac{c}{\sqrt{\eps_0}}
        \biggr)^{n+1} \biggnorm{\sum_{k=1}^K r_k x_k }_{L^2(\Om,E)}\\
& = &  \frac{1}{\eps_0- \eps}\biggnorm{\sum_{k=1}^K r_k  x_k }_{L^2(\Om,E)},
\end{eqnarray*}
where we used that $\eps_0 = \sfrac{1}{4c^2}$. This proves the claim. Now the
result is obtained via \cite[Proposition 2.8]{Weis:FM}.
\end{proof}

As an application of Theorem~\ref{thm:gammaDatkoPazyResolvent} we have the following
bounded perturbation result for the existence of a solution for the perturbed
problem.

\begin{theorem}\label{thm:admiss-perturb0}

Let $P\in \cB(E)$ and $B\in \cB(H,E)$. If $\SCP{A}{B}$ has a solution, then
$\SCP{A{+}P}{B}$ has a solution as well.
\end{theorem}
\begin{proof}
For $\omega\in \RR$ denote $A_{\omega} = A - \omega$ and
$T_{\omega}(\cdot) := e^{-\omega \cdot}T(\cdot)$.
It follows from \cite[Proposition 4.5]{vanNeervenWeis:asymptotic-scp} that for
all $\omega>\omega_0(A)$, $T_{\omega}(\cdot) B\in \ga(\RR_+, H, E)$.
From \cite[Corollary 2.17]{KunstmannWeis:Levico} it follows that for all $\omega>\omega_0(A)+1$,
\[ \RB{ R(\la, A_{\omega})\suchthat  \Re{\la} \geq 0} \leq \frac{c}{\omega-\omega_0(A)-1},
\]
where $c$ is a constant depending only on $(T(t))_{t\geq 0}$. Choose $\omega_1>\omega_0(A)+1$ so large that
$\frac{c}{\omega_1-\omega_0(A)-1}\|P\|<1$.
By \cite[Lemma 5.1]{vanNeervenWeis:asymptotic-scp}, $R(i\cdot,A_{\om_1})B\in
\ga(\RR_+,H,E)$.

Denote by $(S(t))_{t\ge 0}$ the semigroup generated by $A{+}P$
(cf. \cite[Section III.1]{EngelNagel} or \cite[Chapter III]{Pazy})
and let $S_{\omega_1}(t) := e^{-\omega_1t} S(t)$, $t\ge 0$. Since
\[
    \RB{ R(is, A_{\omega_1})P\suchthat s\in  \RR }
\le \RB{ R(is, A_{\omega_1})\suchthat s\in  \RR } \, \norm{P} =: C < 1,
\]
it follows from $i\RR \subseteq \varrho(A_{\omega_1})$ that
$i\RR \subseteq \varrho(A_{\omega_1}+P)$ and
\[
R(is, A_{\omega_1}{+}P)B = \sum_{n=0}^\infty \bigl( R(is, A_{\omega_1})P \bigr)^n R(is, A_{\omega_1})B
=: R_{A,P,\omega_1}(s) R(is, A_{\omega_1})B.
\]
Moreover, as in Theorem~\ref{thm:gammaDatkoPazyResolvent}, and using the fact that
$C<1$, $\{ R_{A,P,\omega_1}(s) \suchthat s\in \RR\}$ is $R$--bounded with constant $\tfrac{1}{1-C}$.
From \cite[Proposition 4.11]{KaltonWeis:square-function-est}
we deduce that
\[
    \norm{ R(i\cdot, A_{\omega_1}{+}P)B }_{\ga(\RR, H, E)}
\le \tfrac{1}{1-C} \norm{ R(i\cdot, A_{\omega_1})B }_{\ga(\RR, H, E)}.
\]
Now \cite[Lemma 5.1]{vanNeervenWeis:asymptotic-scp} shows that $S_{\omega_1}(\cdot) B\in \ga(\RR_+, H, E)$.
It follows from the right ideal property that for all $t>0$,
\[\|S(\cdot) B\|_{\ga(0,t, H, E)} \leq e^{t\omega_1} \|S_{\omega_1}(\cdot) B\|_{\ga(0,t, H, E)}\]
and the result can be obtained via \cite[Theorem 7.1]{vanNeervenWeis:stoch-int}.
\end{proof}

Concerning existence and uniqueness of
invariant measures we obtain:

\begin{theorem}\label{thm:admiss-perturb}

Assume that $s(A) <0$ and that $\{ R(is, A)\suchthat
s\in \RR\}$ is $R$--bounded. Let $B\in \cB(H,E)$ such that
$\SCP{A}{B}$ admits an invariant measure.
Then there exists a $\delta >0$ such that for all $P \in \cB(E)$ with
$\norm{P} < \delta$, $\SCP{A{+}P}{B}$ admits a unique invariant measure.
\end{theorem}
\begin{proof}

Let $\delta>0$ such that $\RB{ R(is, A)\suchthat s\in \RR } \le
\sfrac{1}{\delta}$. Then, if $\norm{P}< \delta$,
\[
    \RB{ R(is, A)P\suchthat s\in  \RR }
\le \RB{ R(is, A)\suchthat s\in  \RR } \norm{P} =: C < 1.
\]

As in Theorem \ref{thm:admiss-perturb0} it can be deduced that

\[
    \norm{ R(i\cdot, A{+}P)B }_{\ga(\RR, H, E)}
\le \tfrac{1}{1-C} \norm{ R(i\cdot, A)B }_{\ga(\RR, H, E)}.
\]
The existence of an invariant measure now follows from \cite[Proposition 4.4
and Lemma 5.1]{vanNeervenWeis:asymptotic-scp}.

By \cite[Corollary 4.3]{vanNeervenWeis:asymptotic-scp}, for uniqueness it
suffices to  note that $R(\la, A+P)$ is uniformly bounded
for $\Re{\la} >0$. \end{proof}

In particular, the $R$-boundedness of $\{R(is, A)\suchthat
s\in \RR\}$ implies that an invariant measure for
$\SCP{A}{B}$, if one exists, is unique.
On the other hand, 
if $i\RR\subseteq\varrho(A)$ but $\{R(is, A)\suchthat
s\in \RR\}$ fails to be $R$-bounded,
then Theorem~\ref{thm:gammaDatkoPazyResolvent}
shows that there exists a rank one operator $B'\in\cB(H,E)$ such that the problem
$\SCP{A}{B'}$ fails to have an invariant measure.
As a result we obtain that if $\SCP{A}{B}$ fails to have a unique
invariant measure, then there exists a rank one operator $B'\in\cB(H,E)$ such that the problem $\SCP{A}{B'}$ fails to have an invariant measure.
A related result can be found in \cite{GN}.

\begin{proof}[{Proof of Theorems {\ref{thm:DatkoPazy}} and {\ref{thm:DatkoPazy2}}}]

If $T(\cdot)x \in \ga(\RR_+,E)$ for all $x\in E$, then by
Theorem~\ref{thm:gammaDatkoPazyResolvent} $s(A) <0$ and
$\{ R(is, A)\suchthat s\in \RR\}$ is $R$--bounded. Thus,

Theorem~\ref{thm:admiss-perturb} applies to the bounded perturbation
$P = \delta \cdot I_E$.
\end{proof}

\def\SUBMITTED{Submitted}
\def\TOAPPEAR{To appear}
\def\PREPARATION{In preparation}

\def\cprime{$'$}
\providecommand{\bysame}{\leavevmode\hbox to3em{\hrulefill}\thinspace}

\contact
\end{document}